\documentclass[12pt]{article}
\usepackage{amsmath,amssymb,theorem}
\usepackage{graphicx, enumerate}
\usepackage{chngcntr}
\usepackage{epstopdf}
\usepackage{setspace}

\counterwithin{figure}{section}
\numberwithin{figure}{section}
\counterwithin{table}{section}
\usepackage{hyperref}
\hypersetup{
  colorlinks,
  citecolor=black,
  filecolor=black,
  linkcolor=black,
  urlcolor=black
}

\newtheorem{thm}{Theorem}[section]
\newtheorem{conj}[thm]{Conjecture}

\newtheorem{lem}[thm]{Lemma}

\theorembodyfont{\rmfamily}

\def\pf{\bigskip\noindent {\bf Proof.}~~}

\def\proofsquare{\bigskip\hfill $\blacksquare$}

\newcounter{SectionStep}[section]

\newcounter{SideStep}[SectionStep]

\begin{document}
\title{Graphs with no $K_9^=$ minor are 10-colorable}
\author{Martin Rolek\thanks{E-mail address: msrolek@wm.edu.} \\
Department of Mathematics\\
College of William \& Mary\\
Williamsburg, VA 23185
}
\date{March 26, 2018}

\maketitle

\begin{abstract}
Hadwiger's conjecture claims that any graph with no $K_t$ minor is $(t - 1)$-colorable.
This has been proved for $t \le 6$, but remains open for $t \ge 7$.
As a variant of this conjecture, graphs with no $K_t^=$ minor have been considered, where $K_t^=$ denotes the complete graph with two edges removed.
It has been shown that graphs with no $K_t^=$ minor are $(2t - 8)$-colorable for $t \in \{7, 8\}$~\cite{Jakobsen1971,Rolek2017}.
In this paper, we extend this result to the case $t = 9$ and show that graphs with no $K_9^=$ minor are $10$-colorable.
\end{abstract}


\section{Introduction}\label{sec:Intro}

All graphs considered are simple and finite.
We use $V(G)$, $|G|$, $E(G)$, $e(G)$, $\delta(G)$, $\alpha(G)$ and $\chi(G)$ to denote the vertex set, number of vertices, edge set, number of edges, minimum degree, independence number, and chromatic number of a graph $G$, respectively.
Given $S \subseteq V(G)$, we denote by $G[S]$ the subgraph of $G$ induced by $S$, and by $G - S$ the subgraph $G[V(G) \setminus S]$ of $G$.
If $S = \{y\}$, we simply write $G - y$ in the latter case.
$G$ contains $H$ as a minor, denoted by $G \ge H$, if $H$ can be obtained from $G$ by some sequence of vertex deletions, edge deletions, and edge contractions.

The famous Hadwiger's conjecture~\cite{Hadwiger1943} claims that every graph with no $K_t$ minor is $(t - 1)$-colorable for $t \ge 1$.
The conjecture has been shown true for $t \le 6$.
The cases $t \le 3$ are easy to show.
The case $t = 4$ was shown by both Hadwiger~\cite{Hadwiger1943} and Dirac~\cite{Dirac1952}, and a short alternative proof has also been given by Woodall~\cite{Woodall1992}.
Wagner~\cite{Wagner1937} proved that the case $t = 5$ is equivalent to the Four Color Theorem, and most recently Robertson, Seymour, and Thomas~\cite{Robertson1993} showed the same for the case $t = 6$ in 1993.
For $t \ge 7$, the conjecture remains open, though there have been some partial results for small values of $t$.
Graphs with no $K_t$ minor for $t \in \{7, 8\}$ were shown to be $(2t - 6)$-colorable by Albar and Gon\c calves~\cite{Albar2017}.
The present author and Song~\cite{Rolek2017} extended this result to the case $t = 9$, and provided an alternative, computer-free proof for the cases $t \in \{7, 8\}$.
In fact, the following conjecture was posed and subsequent theorem was proved in~\cite{Rolek2017}.

\begin{conj}\label{conj:ExtremalColor}
(Rolek and Song~\cite{Rolek2017})
For every $t \ge 1$, every graph on $n$ vertices with at least $(t - 2)n - \binom{t - 1}{2} + 1$ edges either has a $K_t$ minor or is $(t - 1)$-colorable.
\end{conj}

\begin{thm}\label{thm:KtColoring}
(Rolek and Song~\cite{Rolek2017})
For $t \ge 6$, if Conjecture~\ref{conj:ExtremalColor} is true, then every graph with no $K_t$ minor is $(2t - 6)$-colorable.
\end{thm}


As a weaker variant of Hadwiger's conjecture, in this paper we will investigate graphs with no $K_t^=$ minor, where $K_t^=$ denotes the graph obtained from $K_t$ by deleting two edges.
Note that there are two nonisomorphic graphs $K_t^=$, depending on whether the removed edges share a common end or not.
Let $\mathcal{K}_t^=$ denote the family consisting of the two nonisomorphic graphs $K_t^=$.
Throughout this paper, we will use the following conventions.
We say a graph $G$ has no $K_t^=$ minor if $G$ does not contain $K$ as a minor for any $K \in \mathcal{K}_t^=$, and we say that $G$ has a $K_t^=$ minor if $G$ contains $K$ as a minor for some $K \in \mathcal{K}_t^=$.
Jakobsen~\cite{Jakobsen1971} proved the extremal function for $K_7^=$ minors, an easy consequence of which is that graphs with no $K_7^=$ minor are $6$-colorable.
We note that this result is best possible, as evidenced by the $6$-chromatic, $K_7^=$ minor free graph $K_6$.
The present author and Song~\cite{Rolek2017} showed that graphs with no $K_8^=$ minor are $8$-colorable.
The main result of this paper is the following next case.

\begin{thm}\label{thm:K9=Coloring}
Every graph with no $K_9^=$ minor is $10$-colorable.
\end{thm}

To prove Theorem~\ref{thm:K9=Coloring}, we will need the following extremal function for $K_9^=$ minors proved by the present author in~\cite{Rolek2018a}.

\begin{thm}\label{thm:K9=Extremal}
(Rolek~\cite{Rolek2018a})
If $G$ is a graph with $|G| \ge 8$ and at least $6|G| - 20$ edges, then either $G \ge K_9^=$ or $G$ is a $(K_8, K_{2, 2, 2, 2, 2}, 5)$-cockade.
\end{thm}

A $(K_8, K_{2, 2, 2, 2, 2}, 5)$-cockade is a graph built up from disjoint copies of $K_8$ and $K_{2, 2, 2, 2, 2}$ by identifying cliques of size 5.
It is easy to see that any $(K_8, K_{2, 2, 2, 2, 2}, 5)$-cockade is $8$-colorable.

We will additionally need the following, which is an abridged version of a key lemma used in the proof of Theorem~\ref{thm:K9=Extremal} (see Lemma 2.9 in~\cite{Rolek2018a}).
We note that Lemma~\ref{lem:K7=UK1Computer} has been proved by a computer search.

\begin{lem}\label{lem:K7=UK1Computer}
(Rolek~\cite{Rolek2018a})
If $G$ is a graph with $|G| = 11$ and $\delta(G) \ge 6$, then there exists $y \in V(G)$ such that $G - y \ge K_7^=$.
\end{lem}

In our proof of Theorem~\ref{thm:K9=Coloring}, we will investigate a minimum counterexample, chosen among all such graphs to be contraction-critical.
A graph $G$ is $k$-contraction-critical if $\chi(G) = k$, and any proper minor of $G$ is $(k - 1)$-colorable.
Hadwiger's conjecture is equivalent to the claim that the only $k$-contraction-critical graph is the complete graph $K_k$.
Non-complete contraction-critical graphs were first studied by Dirac~\cite{Dirac1960}, and the following is a useful consequence of his initial work on the subject.

\begin{lem}\label{lem:alphaNx}
If $G$ is $k$-contraction-critical, then for any $x \in V(G)$, $\alpha(G[N(x)]) \le d(x) - k + 2$.
\end{lem}

It was also shown by Dirac~\cite{Dirac1960} that $5$-contraction-critical graphs are $5$-connected.
This was extended by Mader~\cite{Mader1968a}, who showed both that $6$-contraction-critical graphs are $6$-connected, and the following deep result.

\begin{thm}\label{thm:MaderConn}
(Mader~\cite{Mader1968a})
For $k \ge 7$, any $k$-contraction-critical graph is $7$-connected.
\end{thm}

It seems difficult to improve Theorem~\ref{thm:MaderConn} for small values of $k$.
Some improvements for larger $k$ have been found, e.g.~\cite{Kawarabayashi2007,Kawarabayashi2013}, with the best current result by Chen, Hu, and Song~\cite{Chen2018+} who showed that any $k$-contraction-critical graph is $\lceil k/6 \rceil$-connected.
In proving Theorem~\ref{thm:K9=Coloring}, we will consider $11$-contraction-critical graphs, and it follows that these graphs will be $7$-connected.
In order to prove Theorem~\ref{thm:K9=Coloring} using connectivity arguments alone, $7$-connectivity would not be sufficient.
Instead, we use the following method introduced by the present author and Song in~\cite{Rolek2017}, which can connect specified nonadjacent vertices in a neighborhood via Kempe chains, or color alternating paths.

\begin{lem}\label{lem:KempeChains}
(Rolek and Song~\cite{Rolek2017})
Let $G$ be a $k$-contraction-critical graph.
Let $x \in V(G)$ be a vertex of degree $k + s$ with $\alpha(G[N(x)]) = s + 2$, and let $S \subset N(x)$ be an independent set with $|S| = s + 2$, where $k \ge 4$ and $s \ge 0$ are integers.
Let $M$ be a set of missing edges of $G[N(x) \setminus S]$.
Then there exists a collection $\{P_{uv} : uv \in M \}$ of paths in $G$ such that for each $uv \in M$, $P_{uv}$ has ends $\{u, v\}$ and all its internal vertices in $G - N[x]$.
Moreover, if the vertices $u$, $v$, $w$, and $z$ are distinct with $uv, wz \in M$, then the paths $P_{uv}$ and $P_{wz}$ are vertex-disjoint.
\end{lem}

\section{Proofs}

In this section, we first prove two useful lemmas.
We then conclude the section with the proof of Theorem~\ref{thm:K9=Coloring}.

\begin{lem}\label{lem:2K7}
Let $G$ be a $7$-connected graph with $|G| \ge 9$.
If $G$ has two different $K_7$ subgraphs, then $G \ge K_9^=$.
\end{lem}

\pf
Let $U_1, U_2 \subseteq V(G)$ such that $G[U_1]$ and $G[U_2]$ are isomorphic to $K_7$, and $U_1 \ne U_2$.
If $|U_1 \cap U_2| = 6$, let $x \in V(G) \setminus (U_1 \cup U_2)$.
Then there exist seven paths $P_1, \dots, P_7$, with one end $x$ and the other end in $U_1 \cup U_2$, which are disjoint except for their common end.
By contracting each of these paths to a single edge, we obtain a $K_9^=$ minor in $G$.
Thus we may assume $|U_1 \cap U_2| \le 5$.
Now there exist seven paths $P_1, \dots, P_7$ with one end in $U_1$ and the other end in $U_2$.
Possibly, some $P_i$ consist of only a single vertex, and we may assume any such path has $i \in \{1, \dots, 5\}$.
By contracting each of $P_1, \dots, P_5$ to a single vertex and $P_6, P_7$ to a single edge, we again obtain a $K_9^=$ minor in $G$.
\proofsquare

\begin{lem}\label{lem:K7=UK1Color} Suppose $G$ is an $11$-contraction-critical graph.
If there exists a vertex $x \in V(G)$ and $y \in N(x)$ such that $d(x) = 11$ and $G[N(x)] - y \ge K_7^=$, then $G \ge K_9^=$.
\end{lem}

\pf
Say $N(x) = \{y, u_0, u_1, \dots, u_9\}$.
If $G[N(x)]$ contains a $K_7$ subgraph, then $G[N[x]]$ contains a $K_8$ subgraph, and in particular $G$ contains two $K_7$ subgraphs.
Since $G$ is $11$-contraction-critical, it is $7$-connected by Theorem~\ref{thm:MaderConn}.
Therefore $G \ge K_9^=$ by Lemma~\ref{lem:2K7}.
Thus we may assume $G[N(x)]$ does not contain $K_7$ as a subgraph.
In particular, we may assume that $G[\{u_0, \dots, u_9\}]$ contains two disjoint missing edges, say $u_0 u_1, u_2 u_3 \notin E(G)$.
Note that by Lemma~\ref{lem:alphaNx}, $\alpha(G[N(x)]) = 2$.
By Lemma~\ref{lem:KempeChains} applied to $N(x)$ a first time with $S = \{u_0, u_1\}$ and $M = \{ yu_i \notin E(G) : 2 \le i \le 9 \}$, and a second time with $S = \{u_2, u_3\}$ and $M = \{ yu_i \notin E(G) : 0 \le i \le 1 \}$, we obtain paths $P_0, P_1, \dots, P_9$ such that for $0 \le i \le 9$ the path $P_i$ has ends $y$ and $u_i$, and $P_i$ consists of only the edge $yu_i$ if $yu_i \in E(G)$.
Note that the paths $P_i$ are not necessarily internally disjoint, but no path $P_i$ has an internal vertex in $N(x)$.
Now by contracting $G[N(x)] - y$ to obtain $K_7^=$, contracting each path $P_i$ to a single edge by contracting onto $y$, and including $x$, we see that $G \ge K_9^=$.
\proofsquare

We are now ready to prove our main result.

\bigskip\noindent{\bf Proof of Theorem~\ref{thm:K9=Coloring}.}
Suppose $G$ is a graph with no $K_9^=$ minor such that $\chi(G) \ge 11$.
We may assume that $G$ is chosen to be contraction-critical.
By Theorem~\ref{thm:K9=Extremal}, $e(G) \le 6|G| - 20$, and it follows that $\delta(G) \le 11$.
Thus $\chi(G) \le 12$.
If $\chi(G) = 12$, then $\delta(G) = 11$, and it follows from Lemma~\ref{lem:alphaNx} that $G[N(x)]$ is isomorphic to $K_{11}$ for any $x \in V(G)$ with $d(x) = 11$, a contradiction.
Therefore, $G$ is $11$-contraction-critical and $\delta(G) \ge 10$.
If $\delta(G) = 10$, then from Lemma~\ref{lem:alphaNx} we see $G[N(x)]$ is isomorphic to $K_{10}$ for any any $x \in V(G)$ with $d(x) = 10$, again a contradiction.
Therefore $\delta(G) = 11$, and once more from Lemma~\ref{lem:alphaNx}, it follows that $\alpha(G[N(x)]) = 2$ for any $x \in V(G)$ with $d(x) = 11$.

Let $x \in V(G)$ such that $d(x) = 11$.
We claim that $G[N(x)]$ contains $K_6$ as a subgraph.
So suppose not.
Since $\alpha(G[N(x)]) = 2$, we must have $\delta(G[N(x)]) \ge 5$.
If $\delta(G[N(x)]) \ge 6$, then it follows from Lemma~\ref{lem:K7=UK1Computer} that there exists $y \in N(x)$ such that $G[N(x)] - y \ge K_7^=$, and so $G \ge K_9^=$ by Lemma~\ref{lem:K7=UK1Color}, a contradiction.
Therefore $\delta(G[N(x)]) = 5$.
Now let $y_1 \in N(x)$ such that $y_1$ is adjacent to five vertices of $N(x)$.
Say $y_1$ is adjacent to $y_2, y_3, \dots, y_6 \in N(x)$, and $N(x) \setminus \{y_1, \dots, y_6\} = \{z_1, \dots, z_5\} = Z$.
Since $\alpha(G[N(x)]) = 2$, $G[Z]$ is isomorphic to $K_5$.
Since $G[N(x)]$ does not contain $K_6$ as a subgraph, we may assume $y_5 y_6 \notin E(G)$.
Then, since $\alpha(G[N(x)]) = 2$, all other vertices of $N(x)$ are adjacent to at least one of $y_5$ or $y_6$.
Suppose $G[\{y_2, y_3, y_4\}]$ is not isomorphic to $K_3$, say $y_2 y_3 \notin E(G)$.
Then, similarly, all other vertices of $N(x)$ are adjacent to at least one of $y_2$ or $y_3$.
In particular, at least three vertices of $Z$ are adjacent to $y_2$, say.
But then by contracting $\{y_1, y_5, y_6\}$ to a single vertex, we see that $G[N(x)] - \{y_3, y_4\} \ge K_7^=$.
Then by Lemma~\ref{lem:K7=UK1Color}, $G \ge K_9^=$, a contradiction.
Thus $G[\{y_2, y_3, y_4\}]$ is isomorphic to $K_3$.

By symmetry, we may assume that $y_5$ has more neighbors in $\{y_2, y_3, y_4\}$ than $y_6$.
Then $y_5$ is adjacent to at least two vertices of $\{y_2, y_3, y_4\}$.
We may further assume by symmetry that if $y_5$ and $y_6$ have the same number of neighbors in $\{y_2, y_3, y_4\}$, then $y_6$ has more neighbors in $Z$.
In either case, $y_6$ is adjacent to at least three vertices of $Z$.
Say $y_6 z_1, y_6 z_2, y_6 z_3 \in E(G)$.
Now, $G[\{y_1, \dots, y_5\}]$ is isomorphic to either $K_5^-$ or $K_5$.
If the former, say $y_2 y_5 \notin E(G)$.
Then, similar to the above, $G[N(x)] - \{y_3, y_4\} \ge K_7^=$, as seen by contracting $\{y_1, y_2, y_5\}$ to a single vertex, and so we obtain a contradiction by Lemma~\ref{lem:K7=UK1Color} again.
Hence we may assume $G[\{y_1, \dots, y_5\}]$ is isomorphic to $K_5$.
Since $G[N(x)]$ does not contain $K_6$ as a subgraph, we may assume $y_6 z_5 \notin E(G)$, and so $y_5 z_5 \in E(G)$ since $\alpha(G[N(x)]) = 2$.
Since $\delta(G[N(x)]) = 5$, we may assume by symmetry that $y_6$ is adjacent to either $y_4$ or $z_4$.
By Lemma~\ref{lem:KempeChains} applied with $S = \{y_5, y_6\}$ and $M = \{y_1 z_1, y_2 z_2, y_3 z_3, y_4 z_4, y_1 z_2\}$, we obtain paths $P_1, \dots, P_5$ such that $P_i$ has ends $y_i, z_i$ for $i \in \{1, \dots, 4\}$, and $P_5$ has ends $y_1, z_2$.
Note that the paths $P_3$ and $P_4$ are each disjoint from all other paths, but $P_1$, $P_2$, and $P_5$ are not necessarily internally disjoint.
Let $P_5^*$ be any subpath of $P_5$ with one end in $V(P_1)$, the other end in $V(P_2)$, and no internal vertices in $V(P_1) \cup V(P_2)$.
We now contract the edge $y_5 z_5$ and paths $P_3$ and $P_4$ each to a single vertex, contract $P_5^*$ to a single edge, contract $P_1$ to a single edge by contracting onto $y_1$, and contract $P_2$ to a single edge by contracting onto $z_2$.
Then, along with $y_6$ and $x$, we see that $G \ge K_9^=$, a contradiction.
This proves the claim that $G[N(x)]$ must contain a $K_6$ subgraph, and therefore $G[N[x]]$ contains a $K_7$ subgraph.

Now since $e(G) \le 6|G| - 20$, it follows that there are at least 40 vertices of degree $11$ in $G$.
In particular, there must be two vertices of degree $11$ which are not adjacent, and so $G$ contains two different $K_7$ subgraphs.
Since $G$ is $7$-connected by Theorem~\ref{thm:MaderConn}, it follows from Lemma~\ref{lem:2K7} that $G \ge K_9^=$, a contradiction.
This contradiction completes the proof.
\proofsquare

\end{document}